\documentclass{amsart}  
\usepackage{epsf,epsfig,amssymb}  
  
\newtheorem{theorem}{Theorem}  
  
\newtheorem{lemma}[theorem]{Lemma}  
\newtheorem{corollary}[theorem]{Corollary}  
  
\theoremstyle{remark}  
\newtheorem*{remark}{Remark}

\newcommand{\Z}{{\bf Z}}   
   
\newcommand{\R}{{\bf R}}   
   
\newcommand{\M}{{\mathcal M}}   
\newcommand{\p}[1]{{\mathcal P}_{#1}}  
\newcommand{\A}[2]{{\mathcal A}_{{#1},{#2}}}  
\newcommand{\Aa}[1]{{\mathcal A}_{{#1}}}  
\newcommand{\aL}{{\mathcal A_{*,*}}}  
\newcommand{\aO}{{\mathcal A_{*}}}  

\newcommand{\vpicstdheight}[1]{\kern-10pt\hbox to20pt{$\vcenter{\epsfig{file=#1.eps,width=15pt}}$}\kern 10pt}  
\newcommand{\vpicbigheight}[1]{\kern-10pt\hbox to40pt{$\vcenter{\epsfig{file=#1.eps,width=30pt}}$}\kern 7pt}  
\newcommand{\vpicweirdheight}[1]{\kern-10pt\hbox to60pt{$\vcenter{\epsfig{file=#1.eps,width=45pt}}$}\kern 7pt}  
  
\title{On invariants of Morse knots}

\author[J.Mostovoy]{Jacob Mostovoy}  
\address{Instituto de Matem\'{a}ticas (Unidad Cuernavaca),  
Universidad Nacional Aut\'{o}noma de M\'{e}xico,  
A.P. 273-3 Cuernavaca, Morelos, MEXICO}  
\email{jacob@matcuer.unam.mx}  
  
\author[T.Stanford]{Theodore Stanford}  
\address{Mathematics Department \\  
United States Naval Academy\\  
572C Holloway Road\\  
Annapolis, MD 21402, USA}  
\email{stanford@nadn.navy.mil}  
                                    
\begin{document}  
  
\begin{abstract}  
We define and study Vassiliev invariants for (long) Morse  
knots.  It is shown that there are Vassiliev invariants  
which can distinguish some topologically equivalent Morse  
knots. In particular, there is an invariant of order 3 for  
Morse knots with one maximum that distinguishes two  
different representations of the figure eight knot. We also  
present the results of computer calculations for some  
invariants of low order. It turns out that for Morse knots  
with two maxima there is a $\Z/2$-valued invariant of order  
6 which is not a reduction of any integer-valued invariant.  

\medskip

\noindent {\em Keywords}: Morse knots, Vassiliev invariants.
\end{abstract}  
  
\maketitle  
  
\section*{Introduction}  
In the last few years knot theorists have studied so many species of  
knots (Legendrian, transverse, holonomic, harmonic, Lissajou, virtual  
knots, etc) that perhaps ``low-dimensional zoology'' will develop  
as a subfield in its own right.  As with bears, buffalos or bumble-bees,  
each species of knots has its own particular niche.  
In particular, our interest in Morse knots stems from the fact that  
Morse knots can be viewed as an intermediate class of  
objects between braids and knots and in this regard they  
appear at least in two important contexts in knot theory.  
  
The first appearance of Morse knots (under the name of  
``plat representations'') was in connection with the plat  
closure of braids in the work of J.Birman \cite{Birman}. The  
main result of \cite{Birman} is a Markov-type theorem for  
the plat closure.  Similarly to the classical Markov theorem,  
the allowed moves on braids fall into two categories: moves  
preserving the number of strands and a stabilisation  
move. Morse knots then are equivalence classes of braids  
under the first class of moves; consequently, the theory of  
Morse knots can be regarded as an ``unstable knot theory''.  
It turned out that there exist topologically equivalent  
(i.e.\ isotopic) Morse knots with the same number of  
strands, or, equivalently, with the same number of maxima of  
the ``height function'', which belong to different classes  
with respect to the equivalence relation on Morse  
knots. Examples of such pairs of knots were given by Birman  
\cite{Birman} and Montesinos \cite{Montesinos}.  
  
Another context where Morse knots are relevant is  
Kontsevich's construction of a universal Vassiliev invariant  
for knots \cite{Kontsevich}. This invariant, known as the  
``Kontsevich integral'', is defined in two steps. First one  
chooses an embedding of a given knot into $\R^3$ as a Morse  
knot and defines what is sometimes called the  
``preliminary'' Kontsevich integral. The preliminary  
integral is invariant under {\it Morse} equivalence of  
embeddings and it has to be normalised in a suitable way to  
obtain a genuine knot invariant. It may seem that the only  
information specific to Morse knots that is lost while  
passing to the stabilised (i.e.\ normalised) version of the  
Kontsevich integral is the number of critical points of the  
``height function'' on the knot.  (In fact, on page 145 of  
\cite{Kontsevich} an erroneous claim is made that the  
complete set of invariants of a Morse knot is its  
topological type together with the number of maxima of the  
height function.) However, certain stabilisation is implicit  
in the definition of the algebra of chord diagrams for  
knots. We will see how to refine the relations on chord  
diagrams in order to make the Kontsevich integral sensitive  
to more subtle information specific to Morse knots.  
  
The main purpose of this paper is to show that some  
``unstable'' information about Morse knots can be captured  
by Vassiliev invariants.  In the first section we consider  
the basic properties of Morse knots.  In particular, it is  
shown that the monoid of Morse knots is ``almost  
commutative'' and Morse knots with one maximum are  
classified.  In section 2 we define the Vassiliev invariants  
and give an example which shows that they can distinguish  
isotopic (but not Morse equivalent) knots. We also introduce  
``semi-stable'' invariants; these, however, turn out to be  
``uninteresting'' in the sense that almost all of them are  
invariants of knots up to isotopy. Section 3 contains rough  
estimates of the stability range for the Vassiliev  
invariants of Morse knots. Finally, in section 4 we present  
some results of computer calculations. Here we encounter a  
rather unexpected phenomenon: there is a $\Z/2$-valued  
Vassiliev invariant of Morse knots with two maxima which is  
not a mod 2 reduction of any integer-valued invariant.  
  
We work with  
long knots, which are somewhat easier   
to deal with than compact knots  
in certain situations.   
We do not attempt to answer the question  
whether long and compact knots lead to the same theory.  
  
We assume the reader to be familiar with the basics of Vassiliev knot   
invariants, suitable references are \cite{BN} and \cite{CD}.

\subsection*{Acknowledgments} We are grateful to Sergei Tabachnikov for  
useful comments. The first author thanks Max-Planck-Institut  
f\"ur Mathematik, Bonn for hospitality during the stay at  
which part of the work presented here was done. The second  
author was partially supported by the Naval Academy Research  
Council.

\section{Morse knots}  
Throughout the paper we will work with oriented long  
knots. A long knot $k$ is a smooth non-singular embedding  
$k:\R\to\R^3$ such that the tangent vector to $k$ tends to $(0,0,-1)$  
when  $|t|$ tends to infinity. Thus we adopt the convention that  
all knots ``point downwards''.  
  
We say that $k$ is a {\em Morse knot} if the height function, that is,   
the vertical component of the function $k(t)$ has only a   
finite number of critical points, all of which are non-degenerate.  
Two Morse knots are {\em Morse equivalent} if one can be   
deformed into the other through Morse knots. It is clear that  
Morse equivalent knots are isotopic.  
  
Denote by $\M_n$ the set of Morse equivalence classes of  
knots with $n$ maxima (or, equivalently, $n$ minima) of  
height function and set $\M=\bigcup\M_i$. Abusing the  
terminology, we will refer to elements of $\M$ as to ``Morse  
knots''; this should not lead to confusion. $\M_0$ consists  
of a single element ${\bf [0]}$ which is the class of the  
embedding $t\to (0,0,-t)$. The sets $\M_n$ are formed by  
Morse knots whose isotopy classes have bridge number less  
than or equal to $n+1$.  
  
On Morse knots there is an operation of connected sum: $k\#  
l$ is a knot formed by ``putting $k$ on top of $l$''.  If  
$k\in\M_n$ and $l\in\M_m$ the sum $k\# l$ belongs to  
$\M_{n+m}$, so Morse knots form a graded monoid under  
connected sum. The grading comes from the number of maxima  
and the identity is ${\bf [0]}$.  
  
Denote by ${\bf [k]}$ the ``$k$-hump'' Morse knot as below:  
\[ \epsffile{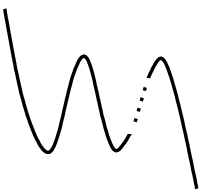} \]  
Notice that ${\bf [m]}\#{\bf [k]} = {\bf [m+k]}$.  
   
The following theorem is a version of a statement in   
\cite{Birman}:  
\begin{theorem}\cite{MoSt}\label{thm:stab}  
For any two isotopic Morse knots $k,l\in\M_n$ there exists  
a non-negative integer $m$ such that $k\#{\bf [m]}$ and $l\#{\bf [m]}$  
are Morse equivalent.  
\end{theorem}  
Thus the theory of Morse knots can be thought of as an unstable   
version of the usual knot theory. Indeed, according to the theorem  
above the direct limit of the sequence  
\[\M_n\stackrel{\#{\bf [1]}}{\to}\M_{n+1}\]  
is precisely the set   
of isotopy classes of knots.  
  
There are other ``unstable knot theories'' which stabilise to  
knots in a similar way.   
For example, closed braids can also be regarded as ``unstable knots''.  
Note, however, that the stabilisation move for closed braids does not  
commute with the unstable moves, that is, with conjugation. In this respect,  
the theory of Legendrian and transverse knots is a better example of an   
unstable knot theory; certain statements and proofs can be transferred from   
the subject of Legendrian and transverse knots to Morse knots   
without major changes.   
In particular, we will apply the technique developed by D.Fuchs  
and S.Tabachnikov in \cite{FT} to describe the semi-stable invariants  
of Morse knots (Theorem \ref{thm:ft} below).

\subsection{Morse knots as closed braids.}  
The stabilisation phenomenon for Morse knots was  
discovered by J.Birman in \cite{Birman} where she proved a   
Markov-type theorem for the   
plat closure. The counterpart of the plat closure  for long knots  
(called ``short-circuit'' closure) is described in   
\cite{MoSt}. The short-circuit closure of a pure braid on  
an odd number of strands is obtained by stretching the top of  
the first strand up and the bottom of the last strand down   
to infinity and joining all strands in turn at the bottom and at   
the top. Every braid in the pure braid group $\p{2n+1}$ closes  
to a knot in $\M_n$. For example, the closure of a trivial braid   
on $2n+1$ strands is the $n$-hump knot ${\bf [n]}$.  
  
There is a Markov theorem for the short-circuit closure:  
\begin{theorem}\cite{MoSt}\label{thm:markov}  
For all $n\geqslant 0$ there exist finitely generated subgroups  
$H^{T}_n, H^{B}_n\subset\p{2n+1}$ such that the fibres of the  
short-circuit map $\p{2n+1}\to\M_n$ are the orbits of the   
simultaneous action on $\p{2n+1}$ by $H^{T}_n$ on the left (that is, on the   
top) and by $H^{B}_n$ on the right.  
\end{theorem}  
  
The stable case of the theorem, i.e.\ the corresponding statement  
for knots up to isotopy is obtained by letting $n$ tend to  
infinity. For an explicit description of $H^{T}_n$ and $H^{B}_n$   
for all $n\geqslant 1$ we refer to \cite{MoSt}. In the last section we list
the generators for  $H^{T}_2$ and $H^{B}_2$; here we will discuss as an   
example the case $n=1$.  
    
\subsection{Morse knots with one maximum.}  
Recall that the pure braid group $\p{3}$ is  isomorphic to  
$F_2\times\Z$ where $F_2$ is the free group on two generators.  
Let $a$ and $b$ be the generators of $F_2$ and $c$ be the   
generator of $\Z$ shown below.  
\[\epsffile{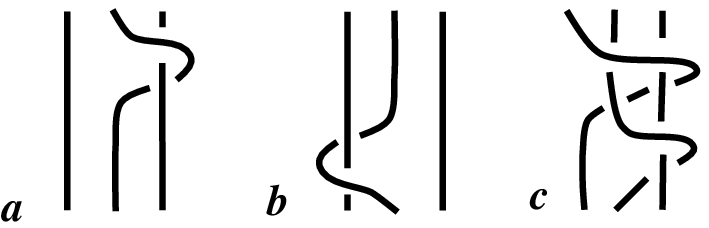}\]  
Then $H^T_1$ is generated by $a$ and $c$, $H^B_1$  
is generated by $b$ and $c$, so the quotient   
$H^T_1\backslash\p{3}/H^B_1=\M_1$ can be easily described:  
\begin{theorem}\label{thm:m1}  
Non-trivial knots in $\M_1$ are in one-to-one correspondence with   
all such reduced words on two letters $a$ and $b$  that  
start with a non-zero power of $b$ and end with a non-zero   
power of $a$.  
\end{theorem}  
By the ``trivial knot'' we mean here the hump ${\bf [1]}$.  
  
Notice that the mirror image of a knot $k$ in $\M_1$ with respect to a   
reflection in any vertical plane can be obtained by replacing $a$ and $b$   
by $a^{-1}$ and  $b^{-1}$ respectively and vice versa in the reduced word  
that represents $k$. This implies the following  
\begin{corollary}  
Any non-trivial knot in $\M_1$ is distinct from its   
mirror image with respect to any vertical plane.  
\end{corollary}  
In particular, the figure eight knot (which is isotopic to its mirror image)  
has at least two distinct representatives in $\M_1$. In the next section we   
will see an example of two Morse knots which are both   
isotopic to the figure eight knot and can be distinguished by Vassiliev   
invariants.  
  
\subsection{Are Morse knots commutative?}  
We do not know if the monoid of Morse knots is commutative.  
The proof of the commutativity of the connected sum operation for  
usual knots consists of making one of the knots very small and running it  
through the other summand. This does not work with Morse knots for the  
following reason: while running one knot through the other one we have turn  
the smaller knot upside down at the maxima and minima of the bigger knot.    
It is not clear if this can be done without creating new critical points.   
Still, we have the following  
\begin{theorem}  
\[\begin{array}{rcl}  
{\rm (a)} & k\# {\bf [1]} = {\bf [1]}\# k &{\rm for\ any\  } k\in\M;\\  
{\rm (b)} & k\# l\# {\bf [1]} = l\# k\# {\bf [1]} & {\rm for\ any\ } k,l\in\M.  
\end{array}\]  
\end{theorem}  
  
Part (a) of the above theorem was proved in  
\cite{MoSt}. Essentially, one has to show that a ``hump''  
can be run through a critical point. The following picture  
illustrates how to do it:  
\[\epsffile{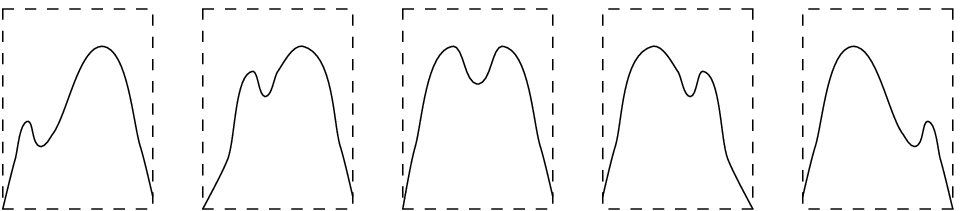}\]  
To verify part (b) notice that the knot $k\# {\bf [1]}$ can be passed through  
a maximum as shown below.  
\[\epsffile{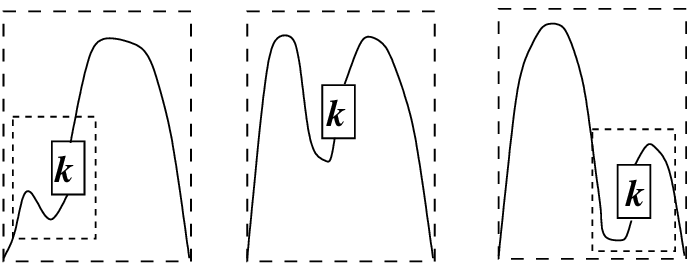}\]  
Let  $k'$ be the Morse knot obtained by passing $k\# {\bf [1]}$ through  
a maximum. Passing  $k'$ through a minimum we again obtain $k\# {\bf [1]}$  
as illustrated by the following picture:   
\[\epsffile{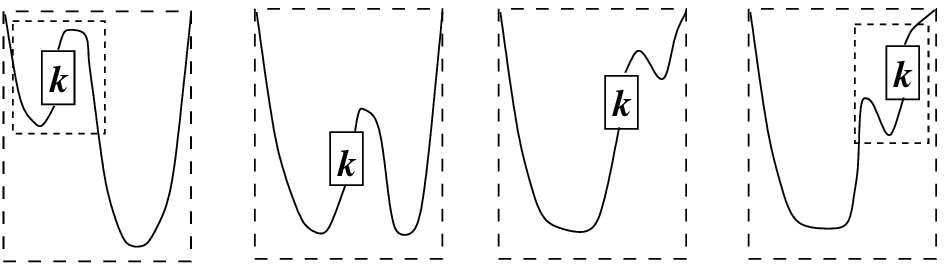}\]  
Passing $k\# {\bf [1]}$ through maxima and minima of $l$ we see that  
$l\# k\# {\bf [1]}$ equals to $k\# {\bf [1]}\# l$ which is the same as  
 $k\# l\# {\bf [1]}$.  
  
\medskip  
  
A possible method of showing that the connected sum is not  
commutative would be to prove the non-commutativity of the  
algebra of Morse chord diagrams defined in the next  
section. Our calculations up to degree $7$ in $\M_2$   
failed to find any such  
phenomenon.  
  
\section{Vassiliev invariants}  
\subsection{Vassiliev invariants of Morse knots.}  
Let $v$ be a function from $\M_n$ to some abelian group  
(which henceforth will be assumed to be the additive group  
of real numbers, unless stated otherwise). It can be  
extended inductively to a function on Morse knots with  
transversal double points by means of the Vassiliev skein  
relation:  
\[ v(\vpicstdheight{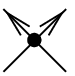}) = v(\vpicstdheight{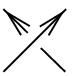}) -v(\vpicstdheight{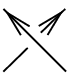}).\]  
The function $v$ is a Vassiliev invariant of type (or order) $k$ if its   
extension vanishes on all knots with more than $k$ double points.  
Notice that we define Vassiliev invariants separately for each   
$\M_n$.  
  
Just like in the usual knot theory one can introduce {\em Morse}  
chord diagrams. A Morse chord diagram with $k$ chords and   
$n$ maxima is a graph as below:  
\[\epsffile{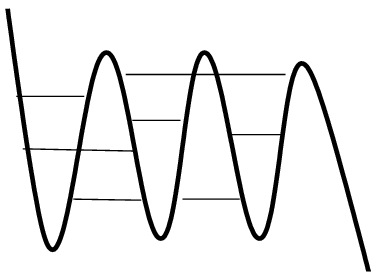} \]  
Here all $k$ chords are horizontal and the Wilson loop is a graph  
of a function on $\R$ with $n$ non-degenerate maxima.  
Let  $\A{k}{n}$ be the \R-vector space generated by all chord   
diagrams with $k$ chords and $n$ maxima modulo several types of  
relations:  
\begin{itemize}  
\item{ homotopy which preserves the number and non-degeneracy of critical   
points and horizontality of chords;}   
\item{ braid-type 4T-relations;}  
\item{ framing independence:  
\[\vpicbigheight{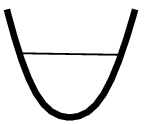} =\vpicbigheight{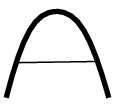}= 0;\]}  
\item{ strand exchange:  
\[\vpicweirdheight{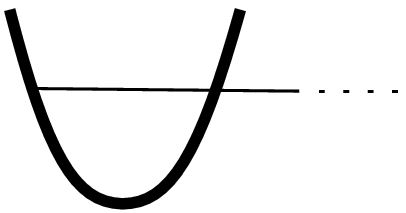}=\vpicweirdheight{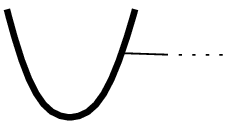}\]  
and  
\[\vpicweirdheight{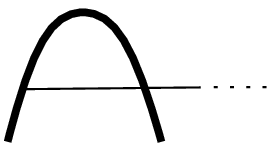}=\vpicweirdheight{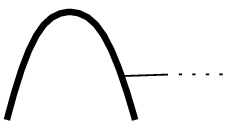} \]}  
\end{itemize}  
Vector spaces $\A{0}{n}$ are one-dimensional; they are generated  
by diagrams without chords and with $n$ maxima.  
  
The connected sum of Morse knots induces a structure of a   
bigraded algebra on  $\aL=\bigoplus_{k,n}\A{k}{n}$.  
\begin{theorem}\label{thm:Kontsevich}  
The Vassiliev invariants of Morse knots form the (graded) dual of the algebra  
$\aL$. In particular, $\A{k}{n}$ is dual to the vector space of the  
invariants of type $k$ modulo invariants of type $k-1$ for  
Morse knots with $n$ maxima.  
\end{theorem}

The tool used to prove this theorem is the universal Vassiliev   
invariant for Morse knots known as the ``Kontsevich integral''.  
The only difference between the Kontsevich integral for Morse  
knots and the Kontsevich integral (without normalisation)  
for knots up to isotopy is that by the ``class of a diagram''   
one should understand its class in the algebra $\aL$ rather than  
$\aO$. The proof is identical to that of the Kontsevich theorem for knots  
up to isotopy apart from that there is no need to check the invariance   
under an insertion of a hump. For details we refer to \cite{BN} or \cite{CD}.  
Note that all terms of the Kontsevich integral  
of a knot in $\M_n$ belong to $\A{*}{n}$. In particular, its first  
term is the chord diagram without chords and with $n$ maxima.  
  
The algebra $\aL$  contains   
some unstable information. For instance, an easy calculation shows that   
the vector space $\A{3}{1}$ is generated by two diagrams  
\[d_1=\vpicbigheight{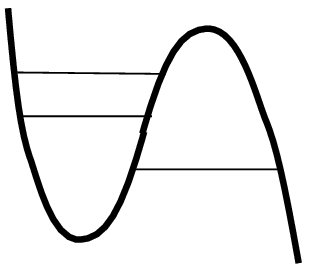}\qquad d_2=\vpicbigheight{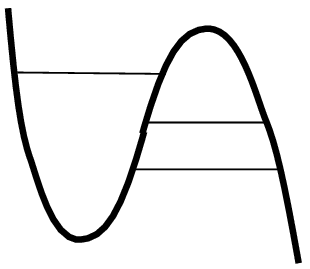}\]  
while there is only one Vassiliev invariant of knots up to   
isotopy in degree 3.  
  
\begin{figure}  
\[\epsffile{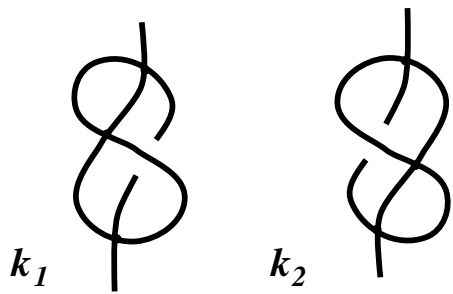}\]  
\caption{Singular knots that are isotopic but not   
Morse equivalent.}  
\label{fig:neq}  
\end{figure}

This implies, for example, that the two singular Morse knots  
$k_1$ and $k_2$ shown on Figure~\ref{fig:neq} are   
distinguished by Vassiliev  
invariants of Morse knots of type 3, as their underlying chord   
diagrams are precisely $d_1$ and $d_2$. On the other hand, their  
isotopy classes differ only by the reversal of orientation.  
It is well-known that all 2-bridge knots are invertible so  
$k_1$ and $k_2$ are in the same isotopy class. Thus Vassiliev invariants  
can distinguish isotopic but not Morse equivalent knots.  
  
In fact, resolving the singularities of both $k_1$ and $k_2$ one obtains 
(up to isotopy) a formal sum of a trefoil and a figure eight knot. 
It is easily checked that a trefoil and its reverse are Morse equivalent  
so, as we have seen in the previous section,   
the figure eight knot has at  
least two distinct representatives in $\M_1$.   
  
\subsection{Semi-stable invariants.}  
The definition of Vassiliev invariants for Morse knots given  
above is not the only possible one. Along with double points  
one can consider yet another type of singularity, namely a  
degenerate critical point where the 3rd derivative of the  
height function does not vanish. The relation for resolving  
degenerate critical points is analoguous to the Vassiliev  
skein relation:  
\[ v(\vpicstdheight{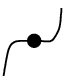})   
= v(\vpicstdheight{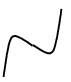}) -v(\vpicstdheight{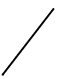}).\]  
We will call the Vassiliev invariants defined in this manner  
{\em semi-stable}.  An example of a semi-stable invariant is  
the number of maxima of a Morse knot, this is an invariant  
of type 1.  It is clear that each semi-stable Vassiliev  
invariant is a Vassiliev invariant. The converse is not true:  
\begin{theorem}\label{thm:ft}   
The algebra of semi-stable Vassiliev invariants of Morse knots is  
generated by Vassiliev invariants of knots up to isotopy and the single   
invariant of type 1, which is the number of maxima.  
\end{theorem}   
The proof is identical to that of Theorems~4.5 and 5.5 of \cite{FT}   
where  similar statements about Legendrian and transverse knots are proved.  
The only modification needed is a substitution of a ``hump'' instead  
of a ``zigzag'' in Legendrian case or a ``double loop'' in the case of   
transverse knots.

\section{Stabilisation of Vassiliev invariants.}  
  
Consider the algebra $\aO$ of usual chord diagrams, i.e.\  
of the chord diagrams for knots up to isotopy. Each Morse chord   
diagram naturally gives rise to a chord diagam of knots  
up to isotopy. This ``forgetful'' map sends    
the relations between Morse chord diagams to relations (or  
trivial identities) in the algebra $\aO$ so   
there is a natural algebra homomorphism   
\[\aL\to\aO\]  
which preserves the grading by the number of chords. The unstable  
Vassiliev invariants of Morse knots correspond to the kernel of this map  
and, as we have seen above, this kernel is non-trivial. It is  
clearly of importance to understand when the Vassiliev invariants  
for Morse knots coincide with those of knots up to isotopy.  
In this section we give a partial answer to   
this question.  
\begin{theorem}\label{thm:Vstab}   
The natural map $\A{k}{n}\to\Aa{k}$ is an epimorphism for   
$k\leqslant 2n$ and an isomorphism for $k\leqslant n$.   
\end{theorem}  
The proof occupies the rest of the section.  
  
We will say that a chord diagram $D$ of knots up to isotopy   
can be presented on $2n+1$ strands if it is the diagram of  
some singular Morse knot with $n$ maxima. The corresponding  
Morse  chord diagram will be called a {\em presentation}  
of $D$ (on $2n+1$ strands).  
  
The first assertion of Theorem~\ref{thm:Vstab}  follows from  
\begin{lemma}\label{lemma:Vstab1}   
Any chord diagram with $2n$ chords can be presented on   
$2n+1$ strands.  
\end{lemma}  
\begin{proof}  
For $n=1$ the statement is   
obvious. Suppose it is established for some $n=p$ and let  
$D$ be a diagram on $2p+2$ chords. Recall that we think of  
knots and diagrams as being ``long'' so the ends of chords of a   
chord diagram are naturally ordered. Consequently, the   
chords themselves can also be ordered according to the order of   
their larger ends.  
Let $c_1$ be the last (i.e.\ the largest) chord of $D$ and   
$c_2$ be the last chord  
of the diagram $D\backslash\{c_1\}$ obtained from $D$ by   
deleting $c_1$. The diagram $D'=D\backslash\{c_1\cup c_2\}$  
can be presented on $2p+1$ strands, so it can also be  
presented on $2p+3$ strands in such a way that no   
chords have ends on any of the last two strands.  
Now it is possible to add the chord $c_2$ to the presentation  
of $D'$ in  such a way that (a) it is horizontal;   
(b) one of its ends is on one of the  first $2p+1$ strands and   
the second end is on the $(2p+2)$nd strand.  
Similarly, one can add the chord $c_1$ to the obtained diagram  
in such a way that its last end is situated on the last strand;  
this gives a presentation of $D$ on $2p+3$ strands.  
\end{proof}  
  
To prove the second assertion of Theorem~\ref{thm:Vstab} it   
suffices to verify the following two lemmas:  
\begin{lemma} \label{lemma:Vstab2}   
If a diagram on $k$ chords can be presented in two different   
ways on $2n+1$ strands, where $n\geqslant k$,    
these  presentations are equivalent modulo the relations   
of strand exchange.   
\end{lemma}  
\begin{lemma}\label{lemma:Vstab3}   
Every 4T and framing independence relation in $\Aa{k}$  
come from relations in $\A{k}{n}$ for $n\geqslant k$;  
\end{lemma}

\begin{proof}[Proof of Lemma~\ref{lemma:Vstab2}]  
Say that a diagram with $n$ chords on $2n+1$ strands is   
{\em arranged} if it has exactly one end of a chord on each of  
the first $2n$ strands. We describe an algorithm of  
connecting any diagram with $n$ chords on $2n+1$ strands  
to an arranged one by relations of strand exchange;  
the lemma follows from the existence of such an algorithm.  
  
Let $N(m)$ be the number of ends of chords on the $m$th strand of  
a given Morse diagram. Suppose the diagram is not arranged.   
  
Find such $m$ that $N(m)=0$ and that $N(i)>0$ for all  
$1\leqslant i<m$. If $N(i)=1$ for all $1\leqslant i<m$ and $m'$ is the   
smallest number such that $m'>m$ and $N(m')\neq 0$ we move one chord end  
from the $m'$th strand to $(m'-1)$st strand by a strand exchange   
relation. (We can do it as $N(m'-1)=0$.) Otherwise, move one  
chord end from $(m-1)$st strand to the $m$th strand.  
  
The above manipulations can be repeated until we get an arranged   
chord diagram. It is a straightforward check that the algorithm  
always terminates on such a diagram.  
\end{proof}  
\begin{proof}[Proof of Lemma~\ref{lemma:Vstab3}]  
A 4T relation can be thought of as a resolution of a ``singular''  
chord diagram whose two chords have a common end. So in order   
to show  
that every 4T relation for diagrams of knots up to isotopy  
comes from a 4T relation for Morse diagrams one can just check  
that every singular chord diagram on $n$ chords can be presented  
on $2n+1$ strands. This can be easily seen using   
Lemma~\ref{lemma:Vstab1}. The statement about the framing  
independence relations is similarly straightforward.   
\end{proof}  
  
\begin{remark}  
The estimates given by   
Theorem~\ref{thm:Vstab} are not sharp.  
Note, however, that nowhere  
in the proof we have used the framing independence relations.  
Also, the statement of the Lemma~\ref{lemma:Vstab3} can be  
easily improved. Namely,  the condition $k\leqslant n$  
can be replaced by $k\leqslant 2n$; the proof is very similar to   
that of Lemma~\ref{lemma:Vstab1}.  
\end{remark}  
  
\section{Computations in $\M_2$.}  
  
We performed some computer calculations of invariants  
of $\M_2$, with the hope of establishing the noncommutativity  
of some Morse knots.  Although we failed in this objective,  
we did discover that there is a $\Z/2$-valued sixth-order  
invariant of $\M_2$ which is not the mod 2 reduction of an  
integer-valued invariant.  Moreover, this invariant can  
distinguish some elements of $\M_2$ from their reverses.  
We describe now how those calculations were carried out.  
  
The basic idea behind our calculations 
is that Vassiliev invariants of Morse knots pull  
back under the short-circuit map to those Vassiliev invariants of pure   
braids which are preserved by the action of the subgroups $H_n^T$ and  
$H_n^B$, see Theorem~\ref{thm:markov}. Thus the problem of computing  
invariants in $\M_n$ is reduced to considering the action of $H_n^T$ and  
$H_n^B$ on pure braid invariants. We refer the reader to \cite{MoWi} for  
a description of the Vassiliev theory for pure braids.   
  
Rather than considering the Vassiliev filtration on knot invariants
we will work with the dual filtration on $\Z$-linear combinations of knots. 
Recall that a knot with $k$ double points can be identified with a formal  
alternating sum of $2^k$ nonsingular knots obtained by resolving all $k$  
double points with the help of the Vassiliev skein relation.  
Consider the abelian group $Q$ generated by all elements of $\M_2$,   
subject to the relations that set all knots with $7$ double points equal to   
zero.  A sixth-order Vassiliev invariant of Morse knots in $\M_2$ 
taking values in an abelian  group $G$ is a 
homomorphism from $Q$ to $G$.  We computed $Q$ to be   
$\Z^{19} \times X$, where $X$ is a finite group whose order is divisible by   
$2$ but not by $3,5,7,11,13,17,19,89$, or $131$. The group $X$ is most  
likely to be $\Z/2$, but since we did all our calculations modulo  
the primes listed, we do not know this for sure.  In any event, there is a   
homomorphism $Q \to \Z/2$ which does not factor through a  homomorphism   
$Q \to \Z$.  
  
The elements of $\M_2$ pull back via the short-circuit map to equivalence   
classes in $P_5$, the equivalence being given by left multiplication by   
$H_2^T$ and right multiplication by $H_2^B$. Denote by $p_{i,j}$ the standard   
generators of $P_5$. The generators of  $H_2^T$ and $H_2^B$ may be chosen as   
below (see \cite{MoSt}):   
   
\smallskip  
  
\begin{tabular}{c|cccccccc}  
$H_2^{T\rule{0pt}{7pt}}$ & $p_{2,3}$  &   
$p_{4,5}$ & $p_{1,2}p_{1,3}$ & $p_{1,4}p_{1,5}$  
& $p_{2,4}p_{2,5}$ & $p_{3,4}p_{3,5}$ & $p_{2,4}p_{3,4}$ & $p_{2,5}p_{3,5}$  
\\[3pt]  
\hline  
$H_2^{B\rule{0pt}{7pt}}$ & $p_{1,2}$ & $p_{3,4}$ & $p_{1,3}p_{2,3}$ & $p_{1,4}p_{2,4}$  
& $p_{1,5}p_{2,5}$ & $p_{1,3}p_{1,4}$ & $p_{2,3}p_{2,4}$ & $p_{3,5}p_{4,5}$  
\\  
\end{tabular}   
  
\medskip  
  
The linear extension of the short-circuit map sends $\Z P_5$ to $\Z$-linear   
combinations of elements of $\M_2$ and knots with $7$ singularities   
are the image of $I^7 \subset \Z P_5$ under this extension.   
(Here $I$ denotes the augmentation ideal of $\Z P_5$.)  We thus take $Q$ to   
be the quotient of $\Z P_5 / I^7$ by the action of $H_2^T$ on the left and   
$H_2^B$ on the right.  
  
Let $q_{i,j} = p_{i,j} - 1$. We call a product of  
$m$ elements from the set  
$\{q_{1,2}, q_{1,3}, \dots q_{4,5}\}$  
{\em a monomial of degree} $m$.  
Let $S$ be the set of monomials of degree less than 7.  
The set $S$ generates $\Z P_5 / I^7$ (as a $\Z$-module),  
subject to relations inherited from those among the  
$p_{i,j}$ in $P_5$. The simplest type of relation among the  
$p_{i,j}$ is a {\em commutation relation}.  For example,  
$p_{1,2} p_{3,4} = p_{3,4} p_{1,2}$.  This translates to  
$q_{1,2} q_{3,4} = q_{3,4} q_{1,2}$.  This relation must  
then be multiplied on the left and on the right by all  
pairs of monomials, the sum of whose degrees is less than $5$,  
to obtain part of the total list of  
relations in a presentation of $\Z P_5 / I^7$.  
There are 10 commutation relations in $P_5$,  
and each generates a similar sublist of the relations  
of $\Z P_5 / I^7$.  
  
Another type of pure braid relation is a {\em three-strand  
relation}, for example $p_{1,2} p_{1,3} p_{2,3} = p_{1,3}  
p_{2,3} p_{1,2}$.  There are 20 such relations in   
$P_5$.  This particular one translates to  
$(q_{1,2}+1)(q_{1,3}+1)(q_{2,3}+1) =  
(q_{1,3}+1)(q_{2,3}+1)(q_{1,2}+1)$, or $q_{1,2} q_{1,3}  
q_{2,3} + q_{1,2} q_{1,3} + q_{1,2} q_{2,3} = q_{1,3}  
q_{2,3} q_{1,2} + q_{1,3} q_{1,2} + q_{2,3} q_{1,2}$.   
The lowest-order terms of such relations are the familiar  
4T braid relations.  
  
The last type of pure braid relation is a {\em four-strand  
relation}.  There are five of these in $P_5$, for example  
$p_{1,3} p_{3,4}^{-1} p_{2,4} p_{3,4} = p_{3,4}^{-1} p_{2,4}  
p_{3,4} p_{1,3}$.  Here things are complicated by the  
presence of inverses of the $p_{i,j}$.  We write  
$p_{3,4}^{-1} = \sum_{i=0}^6 (-1)^i q_{3,4}^i$  
(since we are working modulo $I^7$).    
Our four-strand relation then  
expands to   
$$0 = q_{1,3} q_{2,4} -q_{2,4} q_{1,3}$$  
$$+q_{1,3} q_{2,4} q_{3,4}   
-q_{1,3} q_{3,4} q_{2,4}   
-q_{2,4} q_{3,4} q_{1,3}   
+q_{3,4} q_{2,4} q_{1,3} $$  
$$-q_{1,3} q_{3,4} q_{2,4} q_{3,4}   
+q_{1,3} q_{3,4} q_{3,4} q_{2,4}   
+q_{3,4} q_{2,4} q_{3,4} q_{1,3}   
-q_{3,4} q_{3,4} q_{2,4} q_{1,3} $$  
$$+q_{1,3} q_{3,4} q_{3,4} q_{2,4} q_{3,4}   
-q_{1,3} q_{3,4} q_{3,4} q_{3,4} q_{2,4}   
-q_{3,4} q_{3,4} q_{2,4} q_{3,4} q_{1,3}   
+q_{3,4} q_{3,4} q_{3,4} q_{2,4} q_{1,3} $$  
\begin{flushleft}  
$-q_{1,3} q_{3,4} q_{3,4} q_{3,4} q_{2,4} q_{3,4}   
+q_{1,3} q_{3,4} q_{3,4} q_{3,4} q_{3,4} q_{2,4}$  
\end{flushleft}   
\begin{flushright}  
$+q_{3,4} q_{3,4} q_{3,4} q_{2,4} q_{3,4} q_{1,3}   
-q_{3,4} q_{3,4} q_{3,4} q_{3,4} q_{2,4} q_{1,3} $  
\end{flushright}  
The lowest-order terms of such relations are  
the same as commutation relations.  As with  
the commutation relations, each three-strand and four-strand  
relation must be multiplied on the left and right  
by pairs of monomials, ignoring  
of course those terms of degree greater than $6$.  
  
Then we must consider the relations generated by the action  
of $H_2^T$ and of $H_2^B$.  Multiplying a braid on  
the left by $p_{2,3}$, for example, does not change the  
Morse knot represented.  This translates to any element  
of $S$ which begins with $q_{2,3}$ being equal to $0$.  
Similarly, any element of $S$ which begins with $q_{4,5}$  
or ends with $q_{1,2}$ or $q_{3,4}$ is also $0$.  
We may call such relations {\em topological framing independence  
relations}.  
Multiplying on the left by $p_{1,2} p_{1,3}$  
doesn't change the Morse knot represented, and this implies  
all relations in $Q$ of the form   
$q_{1,2} x + q_{1,3} x + q_{1,2} q_{1,3} x = 0$, where  
$x$ is a monomial of degree less than $5$.  
There are five more such relations from $H_2^T$ and six  
from $H_2^B$. We may call such relations {\em topological strand exchange  
relations}.    
  
We have now described all the relations we need  
for a presentation of $Q$.  
However, the set $S$ has $1,111,111$ elements, and rather than  
attempt to solve over this many unknowns, we consider the set  
$T \subset S$ of all monomials which are generated by the set  
$\{q_{1,5}, q_{2,5}, q_{3,5}, q_{4,5}\}$, and moreover  
which begin with either $q_{1,5}$  
$q_{3,5}$ and end with either $q_{2,5}$ or $q_{4,5}$.  
There are $1365 = 1 + 0 + 4 + 16 + 64 + 256 + 1024$  
(breaking it down by degree) elements of $T$.  
We describe briefly an algorithm for writing any  
element of $S$  
as a linear ($\Z$-linear, as always in this section)  
combination of elements of $T$.  
  
Given a monomial $x$, we first replace it with a  
linear combination of monomials which do not contain  
$q_{1,2}$.  If $q_{1,2}$ occurs at the end of $x$,  
then $x=0$, and we are done.  If $q_{1,2}$ occurs  
in the middle of $x$, then there exists a commutation  
relation or a three-strand relation  
which  
moves $q_{1,2}$ closer to the end of $x$  
at the expense of adding several more linear   
terms, each of  
which is either of higher degree than $x$ or contains  
one fewer occurences of $q_{1,2}$.  For example,  
$q_{1,2}$ moves past $q_{3,4}$ with a commutation  
relation and past $q_{2,3}$ with a three-strand relation.   
We may then proceed  
inductively to eliminate all occurences of $q_{1,2}$.  
  
The next step is to eliminate $q_{1,3}$ and $q_{2,3}$.  
Both of these may be moved toward the end of a  
monomial $x$ by a process similar to that just described.  
Four-strand relations may now be necessary, for example  
to move $q_{1,3}$ past $q_{2,4}$.  
One has to be careful not to reintroduce $q_{1,2}$,  
or at least to raise the degree of the monomial if it  
is introduced, but this is always possible.    
  
If $q_{1,3}$ occurs at the end of $x$, then a  
topological strand exchange will change it to $q_{1,4}$  
(a similar exchange will turn $q_{2,3}$ to $q_{2,4}$),  
introducing an extra term whose degree is one greater  
than that of $x$.  Similarly, $q_{i,4}$ may be moved  
toward the beginning of $x$, and then exchanged for  
$q_{i,5}$.  Now all the $q_{i,j}$ occuring in the  
linear expansion of $x$ so far have $j=5$.  Any monomial  
which begins with $q_{4,5}$ is $0$, and any monomial  
which begins with $q_{2,5}$ may be exchanged for one  
which begins with $q_{3,5}$.  Topological strand exchanges  
may also be employed to ensure that $x$ ends in   
$q_{2,5}$ or $q_{4,5}$.  
  
We took all the relations---commutation, three-strand,  
four-strand, and topological strand exchange and framing  
independence---and applied the above algorithm to write  
each one as a linear combination of elements of $T$.  This  
was the most time-consuming part of the calculation,  
taking several days using compiled C code on a 300 MHz  
processor.  It was not unusual for a single monomial in  
$S$ to require several hundred million iterations  
of the algorithm to resolve it into a linear combination  
of monomials in $T$.   
What we did was overkill, of course, since many   
of the relations that we processed were accounted for  
by being used in the algorithm.  In fact, most of the  
relations came out to be trivial in the end, though  
there still remained several   
hundred thousand nontrivial relations  
in the 1365 variables.    
  
Exact rational  
solution of these equations  
did not seem feasible, so we solved these  
equations modulo various primes.  For every  
prime we tried except $2$, the dimension of the  
solution space was $19$.  For $2$, the dimension  
was $20$.  Thus we obtain $Q = \Z^{19} \times X$,   
as noted above.  The reason  
we know that $\Z^{19}$ is really $\Z^{19}$, and not  
$19$ factors of the form $\Z/P$ for various large  
numbers $P$ divisible by many small primes, is that   
$19 = 1+0+1+1+3+4+9$, and these numbers are the dimensions  
of the spaces of rational Vassiliev invariants (of isotopy).  
Moreover, we verified in our calculations that  
the map  $\A{6}{2}\to\Aa{6}$ described in Section 3  
is an epimorphism.  
  
To investigate what happens when a knot is reversed,  
we generated all the relations of the form  
``$x$ equals its reverse'', where $x \in S$.  
(The reverse of Morse knot is obtained by rotating  
the braid presenting it $180^\circ$ in the plane,  
and this rotation induces an obvious involution  
on the set $S$.)  We processed these relations  
with the same algorithm to obtain relations among  
the elements of $T$.  When we added these relations  
to all the previous  
ones, the dimension of the solution space  
was $19$ for all primes, including $2$, thus indicating  
that the mod $2$ invariant which we had found sometimes  
detects the difference between a knot and its  
reverse.  
  
We note that there are various modifications   
of the algorithm   
described.  It is possible, for example,   
to avoid using the framing independence relation at all.  
One must add monomials that begin with $q_{4,5}$ to the  
set $T$, and instead of pushing $q_{1,2}$ to the end  
of a monomial $x$, as described above, one pushes it  
to the beginning and exchanges it for   
$-q_{1,3}-q_{1,2}q_{1,3}$.  We attempted to run this  
version of the algorithm, but (in degree six) it   
was taking many times as long as the version described  
above, and so we gave up.

\bibliographystyle{amsplain}  
\bibliography{morse}   
\end{document}